\documentclass[a4paper,12pt]{article}
\usepackage[english]{babel}
\usepackage[T2A]{fontenc}
\usepackage[cp1251]{inputenc}
\usepackage{amsthm}
\usepackage[tbtags]{amsmath}
\usepackage{amsfonts,amssymb}
\begin{document}

\newtheorem{theorem}{Theorem}
\newtheorem{lemma}{Lemma}
\newtheorem{proposition}{Proposition}
\newtheorem{Cor}{Corollary}


\begin{center}
{\large\bf On Centrally Essential Rings}
\end{center}

\begin{center}
Oleg Lyubimtsev\footnote{Nizhny Novgorod State University, Nizhny Novgorod, Russia; email: oleg\_lyubimcev@mail.ru .},
Askar Tuganbaev\footnote{National Research University MPEI, Moscow, Russia; Lomonosov Moscow State University, Moscow, Russia; tuganbaev@gmail.com .}
\end{center}

\textbf{Abstract.} A ring $R$ with center $C$ is said to be centrally essential if the module $R_C$ is an essential extension of the module $C_C$. In this paper, we study properties of ideals of centrally essential rings, centrally essential quaternion algebras, and group rings of Hamiltonian groups.

\textbf{Key words.} centrally essential ring, quaternion algebra, group ring

The work of Oleg Lyubimtsev is supported by Ministry of Education and Science of the Russian Federation, project FSWR-2023-0034. The study of Askar Tuganbaev is supported by grant of Russian Science Foundation (=RSF), project 22-11-00052, https://rscf.ru/en/project/22-11-00052.

\textbf{MSC2020 database 16D25, 16R99}

\section{Introduction}\label{sec1}

We consider only associative rings which are not necessarily unital. A ring $R$ is said to be \textsf{centrally essential} if $R$ is either  commutative or for any its non-central element $a$, there exist non-zero central elements $x$ and $y$ with $ax = y$. It is clear that a unital ring $R$ with center $C(R)$ is centrally essential if and only if the module $R_{C(R)}$ is an essential extension of the module $C(R)_{C(R)}$. An unital ring $R$ is centrally essential if and only if for any non-zero element $a\in R$, there exist non-zero central elements $x,y\in R$ with $ax = y$. The last condition does not hold for any non-zero ring with zero multiplication (it is centrally essential). We note that any commutative ring is centrally essential. 

In Section 2, we study properties of ideals of centrally essential rings which are semiprime as rings. We prove that ideals (resp., right ideals) of a centrally essential ring, which are semiprime as rings, are contained in the center of the ring (resp., are reduced commutative ideals); see Theorem 2.2 and Proposition 2.4. As a corollary, we obtain that every minimal non-nilpotent right ideal of a centrally essential ring is contained in the center of the ring.

In Section 3, we describe centrally essential quaternion algebras and centrally essential group rings of Hamiltonian groups (Theorems 3.2 and 3.6). 

In Section 4, we prove commutativity of a ring which is a rational extension of its center (Proposition 4.2) and give examples of non-commutative centrally essential rings with injective right maximal ring of fractions (Example 4.3) and a non-commutative centrally essential semiring of order 5 (Example 4.4). 

We denote by $Q_{cl}(R)$ and $Q_{max}(R)$ the classical ring of fractions and the maximal (right) ring of fractions of the ring $R$, correspondingly. A right ideal $I$ of the ring $R$ is said to be \textsf{essential} if $I\cap J\neq 0$ for all non-zero right ideals $J$ of $R$. We denote by $\text{Ann}_R(S)$ the left annihilator $\{r\in R \,|\, rS = 0\}$ of the set $S$ in the ring $R$. The right annihilator of the set $S$ is similarly defined. The right singular ideal $Z_r(R)$ is the set of all elements whose right annihilators are essential. A ring $R$ is said to be \textsf{right nonsingular} if $Z_r(R) = 0$. A ring is called a (right) \textsf{ring of fractions} of $R$ if it contains $R$ and is contained in $Q_{max}(R)$. It is well known that a ring $R$ has a von Neumann regular maximal ring of fractions if and only if $R$ is right nonsingular; e.g., see \cite[Proposition 8.3]{F67}. A ring is said to be \textsf{right (resp., left) strongly bounded} if every right (resp., left) ideal contains a non-zero ideal; e.g., see \cite{BT88}. A ring is said to be \textsf{strongly bounded} if is right and left bounded. We set $[a,b] = ab - ba$ for any two elements $a,b$ of the ring $R$. 

A ring $A$ is said to be \textsf{right invariant} (resp., \textsf{right quasi-invariant}) if every right ideal (resp., maximal right ideal) of $A$ is an ideal of $A$.

We denote by $Q_8$ the quaternion group of order 8, i.e., the group with two generators $a$, $b$ and defining relations 
$a^4 = 1$, $a^2 = b^2$ and $aba^{-1} = b^{-1}$; e.g., see \cite[Section 4.4]{H59}. We have: 
$$
Q_8 = \{e, a, a^2, b, ab, a^3, a^2b, a^3b\},
$$
the center $Z(Q_8) = \{e, a^2\}$.

\section{Properties of Right Ideals which Are Semiprime as Rings}\label{sec2}

We recall that a module $X$ is called an \textsf{essential extension} of a submodule $Y$ in $X$ if $Y\cap Z\ne 0$ for any non-zero submodule $Z$ in $X$. In this case, $Y$ is called an \textsf{essential submodule} of the module $X$. A right ideal $I$ of the ring $R$ is said to be \textsf{reduced} if $I$ does nor contain non-zero nilpotent elements. 

\textbf{Lemma 2.1.} Let $R$ be a centrally essential ring and let $J$ be a right ideal in $R$ which is semiprime as a ring. Then $J$ is a reduced right ideal.

\textbf{Proof.} Let $0\neq a\in J$ and $a^2 = 0$. Since $R$ is a centrally essential ring,  $0\neq ac = d\in C(R)\cap J$, where $c\in C(R)$. Then $dR$ is an ideal in $J$ and $(dR)^2 = 0$. This is a contradiction.~$\square$ 

\textbf{Theorem 2.2.} Let $R$ be a centrally essential ring and let $I$ be an ideal in $R$ which is semiprime as a ring. Then $I\subseteq C(R)$. In particular, $I$ is a centrally essential ring.

\textbf{Proof.} We prove that $I$ is a commutative ideal. Let $[a, b]\neq 0$ for some $a, b\in I$. It follows from \cite[Lemma 1.1(iv)]{ABP09} that $C(R)\cap I = C(I)$. Then $J = \{c\in C(R)\, | \, ac\in C(I)\}$ is a non-zero ideal in $C(R)$. Let $d\in C(I)$ be with $dJ = 0$. If $ad\neq 0$, then  $(ad)z\in C(I)\backslash \{0\}$ for some $z\in C(R)$. Then $dz\in J$ and $d(dz) = (dz)^2 = 0$. It follows from Lemma 2.1 that $dz = 0$, whence $adz = 0$. This is a contradiction. Therefore, $ad = 0$ and $d\in J$. Then $d^2 = 0$ and $d = 0$. Consequently, $\text{Ann}_{C(I)}(J) = 0$, where $\text{Ann}_{C(I)}(J) = \{d\in C(I) \,|\, dJ = 0\}$. Further, we note that for any $j\in J$, we have $aj = ja\in C(R)$. Then
$$
[a, b]j = (ab - ba)j = a(bj) - b(aj) = ajb - ajb = 0.
$$
Therefore, $[a, b]J = 0$. In addition, $c_1[a, b] = c_2\in C(I)\backslash \{0\}$ for some $c_1\in C(R)$. Then $c_2J = c_1[a, b]J = 0$ and 
$\text{Ann}_{C(I)}(J)\neq 0$. This is a contradiction. Consequently, $[a, b] = 0$ and $I$ is a commutative ideal. Since $I = C(I) = I\cap C(R)$, we have $I\subseteq C(R)$.~$\square$ 

\textbf{Corollary 2.3.} Every minimal non-nilpotent right ideal $J$ of a centrally essential ring $R$ is contained in the center of the ring $R$ and is ideal in $R$.

\textbf{Proof.} In \cite[Theorem 1.7(1)]{LT23}, it is proved that $J$ is a two-sided ideal. Since $J^2\neq 0$, the ideal $J$ is semiprime as ring.~$\square$ 

\textbf{Proposition 2.4.} Let $R$ be a centrally essential ring and let $J$ be a right ideal in $R$ which is semiprime as a ring. Then $J$ is a reduced commutative right ideal.

\textbf{Proof.} It follows from Lemma 2.1 that $J$ is a reduced right ideal. In \cite[Proposition 2.6.]{LT23}, it is proved that any centrally essential ring is strongly bounded. It follows from \cite[Proposition 1]{BT88} that $J$ is an essential extension of some ideal $I$. By Theorem 2.2, the ideal $I$ is contained in the center rings $R$. We note that $\text{Ann}_J(I) = \{ j\in J\, |\, jI = 0\} = 0$, otherwise $\text{Ann}_J(I)\cap I\neq 0$, and $I$ is not a reduced ring. Consequently, for every $j\in J$, there exists an element $i\in I$ such that $0\neq ji\in I$. Then $J$ is a centrally essential semiprime ring which have to be commutative; see \cite[Proposition 2.1]{LT21b}.~$\square$ 

Under an additional condition, the right ideal $J$ from Proposition 2.4 is contained in the center of the ring.

\textbf{Corollary 2.5.} Let $R$ be a centrally essential ring and let $J$ be a right ideal of $R$ which is semiprime, as a ring. If 
$\text{Ann}_R(J) = 0$, then $J\subseteq C(R)$.

\textbf{Proof.} Let $x\in R$ and $a, b\in J$. It follows from Proposition 2.4 that 
$$
b(ax - xa) = b(ax) - b(xa) = (ba)x - (bx)a = a(bx) - a(bx) = 0.
$$
Therefore, $ax - xa = t\in \text{Ann}_R(J)$, where $\text{Ann}_R(J)$ is the right annihilator $J$. In addition, $\text{Ann}_R(J)$ also is the left annihilator. Indeed, if $tJ\neq 0$, then $(tJ)^2 = 0$ which contradicts to the property that $J$ is reduced. Therefore, $t = 0$ and $J\subseteq C(R)$. ~$\square$

\textbf{Proposition 2.6.} If $R$ is a centrally essential ring and $I$ is a minimal non-nilpotent right ideal, then $I$ and $R/I$ are centrally essential rings\footnote{We recall that $I$ is an ideal in $R$, by Corollary 2.3.}.

\textbf{Proof.} By assumption, $I = eR$ for some idempotent $e\in C(R)$; see \cite[Lemma 2.3]{MT18}. Then 
$R = eR\oplus (1 - e)R$ is a ring direct sum. Consequently, $R$ is a centrally essential ring if and only if $I = eR$ and 
$J = (1 - e)R$ are centrally essential rings.~$\square$ 

We note that Proposition 2.6 is not true for a nilpotent minimal ideal. 
Let $Q_8$ be the quaternion group, $\widehat Q_8 = \sum_{g\in Q_8}g$,  $R = \mathbb{Z}_2Q_8$, and let $A$ be the factor ring of the centrally essential group ring $R$ with respect to the least non-zero ideal $I = \{0, \widehat Q_8\}$. In the paper \cite[Example 3.2]{LT23}, it is proved that the ring $A$ is not centrally essential.

\textbf{Corollary 2.7.} Let $R$ be a ring, in which all idempotents are central, and let $I$ be a minimal non-nilpotent right ideal in $R$. A ring $R$ is centrally essential if and only if $I$ and $R/I$ are centrally essential rings.

Classes of rings in which idempotents are central include, for example, invariant rings or $UC$-rings; see \cite[Proposition 1]{T60} and \cite[Lemma 4]{NZ04}, correspondingly. From \cite[Proposition 2]{T60}, we obtain the following corollary.

\textbf{Corollary 2.8.} Any non-nilpotent minimal ideal of an invariant centrally essential ring is a field.

We note that there exist centrally essential rings which are not left or right invariant (quasi-invariant); see \cite[Example 2.4]{LT21a} and \cite[Theorem 1.7(2)]{LT23}.

\textbf{Open Question 2.9.} Is it true that a nilpotent minimal ideal of a centrally essential ring is contained in the center of the ring?

\section{Centrally Essential Quaternion Algebras and Group Rings}\label{sec3}

The following assertion is well known, but it is given with a proof for convenience.

\textbf{Lemma 3.1.} Let $R$ be a ring and let $\{M_\alpha\}_{\alpha\in \mathcal{A}}$, $\{N_\alpha\}_{\alpha\in \mathcal{A}}$ be two sets of right $R$-modules such that for any $\alpha\in \mathcal{A}$, the module $N_\alpha$ is an essential submodule of the module $M_\alpha$. Then the module  $N = \oplus_{\alpha\in \mathcal{A}}N_\alpha$ is an essential submodule of the module $M = \oplus_{\alpha\in \mathcal{A}}M_\alpha$.

\textbf{Proof.} For any $\alpha\in \mathcal{A}$, let $\pi_\alpha$ be the canonical epimorphism from $M$ onto $M_\alpha$. For any $x\in M$, we set
$$
k(x)=|\{\alpha\in \mathcal{A}: \pi_\alpha(x)\not\in N_\alpha\}|.
$$
The following properties of the function $k(x)$ are obvious: \\
$\forall x\in M,\quad k(x)<\infty$;\\
$\forall x\in M,\;t\in R,\quad k(tx)\leq
k(x)$;\\
$\forall x\in M,\quad k(x)=0\Leftrightarrow x\in N$.

Now we prove that $xR\cap N\neq 0$ provided $x\in M\setminus\{0\}$. We assume the contrary. Among elements $x\in M\setminus\{0\}$ with $xR\cap N = 0$, we take an element $x_0$ with least value of the parameter $k(x)$. It is clear that $k(x_0)>0$, since $x_0\not\in N$. Therefore, there exists a subscript
$\overline{\alpha}\in \mathcal{A}$ such that $\pi_{\overline{\alpha}}(x_0)\not\in N_{\overline{\alpha}}$. Since $N_{\overline{\alpha}}$ is essential in $M_{\overline{\alpha}}$, there exists an element $a_0\in R$ such that
$\pi_{\overline{\alpha}}(x_0)a_0\in N_{\overline{\alpha}}\setminus\{0\}$. Then it is clear that $x_0a_0\neq 0$ and $k(a_0x_0) < k(x_0)$. Therefore, it follows from the choice of the element $x_0$ that $x_0a_0R\cap N\neq 0$. Since $a_0R\subseteq R$, we have $x_0R\cap N\neq 0$, this is a contradiction. Consequently, $N$ is an essential submodule in $M$.~$\square$

Let $A$ be a commutative ring, $I = Ann_A(2)$, and let $a,b$ be two invertible elements of the ring $A$. By following \cite{tugan93}, we denote by $R = (a,b,A)$ the algebra over $A$ with basis $\{1,i,j,k\}$ such that the multiplication in $R$ is defined by relations
$$
i^2 = a,\; j^2 = b,\;
 ij = -ji = k,\; ik = -ki = aj,\; kj = -jk = bi.
$$
It is known that the center of the ring $R$ is of the form $C(R) = A\cdot 1+I\cdot i+I\cdot j+I\cdot k$; see \cite[Lemma 2(b)]{tugan93}.

\textbf{Theorem 3.2.} The following conditions are equivalent:

$(a)$ the ring $(a,b,A)$ centrally essential;

$(b)$ $Ann_A(2)$ is an essential ideal of the ring $A$.

\textbf{Proof.} $(b)\Rightarrow (a)$. We set $\mathcal{A}=\{0,1,2,3\}$, $M_0 = N_0 = A\cdot 1$, $M_1 = A\cdot i$, $M_2 = A\cdot j$, 
$M_3 = A\cdot k$, $N_1 = I\cdot i$, $N_2 = I\cdot j$, and $N_3 = I\cdot k$. We have $C(R) = N_0\oplus N_1\oplus N_2\oplus N_3$, and it follows from the condition $(b)$ that $C(R)$ is an essential $A$-submodule of the module $R$; $C(R)$ especially is an essential $C(R)$-submodule of $R$.

$(a)\Rightarrow (b)$. Let $t\in A\setminus \{0\}$. We consider the element $r = 2t\cdot i\in R$. If $2t = 0$, then $t\in At\cap I$, whence $At\cap I\neq 0$.
Otherwise, there exists an element $c\in C(R)$ such that $0\neq cr\in C(R)$. We have $c = c_0\cdot 1 + c_1\cdot i + c_2\cdot j + c_3\cdot k$, where $c_0\in A$, $c_1,c_2,c_3\in I$. Consequently, $cr = 2c_0t \cdot i$, whence $2c_0t\neq 0$. Then $2c_0t\in I$, whence $0\neq 2c_0t\in{}At\cap I$. We have that $At\cap I\neq 0$ in any case, i.e., condition $(b)$ holds.~$\square$

In the paper \cite[Theorem 4]{S76}, it is proved that if $F$ is a field of characteristic 0 and $G$ is any group, then the group ring $FG$ is right nonsingular. In the following proposition, we extend this result to commutative domains of characteristic $0$.

\textbf{Proposition 3.3.} Let $R$ be a commutative domain of characteristic $0$. Then the group ring $RG$ is right nonsingular for any group $G$.

\textbf{Proof.} Let $Q = Q_{cl}(R)$ be the field of fractions of the ring $R$. We note that $Q = Q_{max}(R)$; e.g., see \cite[Exercise 2, P. 167]{P77}. It follows from \cite[3.3]{B69} that $QG\subseteq Q_{max}(RG)$. In addition, the group ring $QG$ is nonsingular; see \cite[Theorem 4]{S76}. Since $QG$ is the ring of fractions of $RG$, we have that $QG$ and $RG$ have the same maximal ring of fractions which is von Neumann regular. Consequently, the group ring $RG$ is
nonsingular.~$\square$ 

Centrally essential group rings over fields were studied in the paper \cite{MT18}, where it was proved that centrally essential group algebras over fields of characteristic $0$ are commutative. From \cite[Proposition 2.8]{MT19b}, we obtain the following corollary.

\textbf{Corollary 3.4.} Let $R$ be a commutative domain of characteristic $0$ and let $G$ be any group. The group ring $RG$ is centrally essential if and only if the ring $RG$ is commutative.~$\square$ 

We also have a similar assertion for non-commutative domains of characteristic $0$.

\textbf{Proposition 3.5.} Let $R$ be a domain of characteristic $0$ and let $G$ be any group. The group ring $RG$ is centrally essential if and only if $RG$ is commutative.

\textbf{Proof.} If $char R = 0$ and $R$ is a domain, then $RG$ does not have nil-ideals; in particular, it is semiprime, see \cite[Theorem 5.1]{ZM75}. It follows from \cite[Proposition 3.3]{MT18} that $RG$ is a commutative ring.~$\square$ 

We recall that a group $G$ is said to be \textsf{Hamiltonian} if $G$ is not Abelian and all subgroups of $G$ are normal. If $H$ is a subgroup of the group $G$, then we denote by $\Delta(G, H)$ the fundamental ideal of the subgroup $H$ in the group ring $RG$.

We recall that $Q_8$ denotes the quaternion group of order 8.

\textbf{Theorem 3.6.} Let $G$ be a Hamiltonian group and let $R$ be a unital ring. The group ring $RG$ is centrally essential if and only if $R$ is a centrally essential ring of characteristic $2$.

\textbf{Proof.} It is well known\footnote{See \cite[12.5.4]{H59}.} that a Hamiltonian group is the direct product $G = Q_8\times A\times B$, where $A$ is an elementary Abelian $2$-group and $B$ is an Abelian group in which every element is of finite order. Then
$$
RG = R(Q_8\times A\times B)\cong R(A\times B)Q_8.
$$
It follows from \cite[Lemma 2.1 and Lemma 2.2]{MT18} that the group ring $R(A\times B)Q_8$ is centrally essential if and only if $RQ_8$ is a centrally essential group ring. Therefore, the study of central essentiality of the ring $RG$ is reduced to the case $G = Q_8$. 

Let $RQ_8$ be a centrally essential ring. Then $R$ is a centrally essential ring; see \cite[Lemma 2.1]{MT18}. We assume that $char R > 2$. Then $RQ_8\cong R(Q_8/Q'_8)\oplus \Delta(Q_8, Q'_8)$ is a direct sum; see \cite[Proposition 3.6.7]{MS02}. Therefore, $RQ_8$ is a centrally essential ring if and only if $R(Q_8/Q'_8)$ and $\Delta(Q_8, Q'_8)$ are centrally essential rings. Since $Q_8/Q'_8$ is an Abelian group, 
$R(Q_8/Q'_8)$ is a centrally essential ring; see \cite[Lemma 2.2]{MT18}. We consider the ring $\Delta(Q_8, Q'_8)$ with basis 
$\{f, af, bf, abf\}$, where $f = 1 - e_{Q'_8} = \frac{1}{2}(1 - a^2)$. We describe the center $C = C(\Delta(Q_8, Q'_8))$ of the ring $\Delta(Q_8, Q'_8)$. Let
$\alpha = \alpha_0f + \alpha_1af + \alpha_2bf + \alpha_3abf\in C$, $\beta = \beta_0f + \beta_1af + \beta_2bf + \beta_3abf\in \Delta(Q_8, Q'_8)$. The condition $\alpha\beta = \beta\alpha$ holds if and only if relations
$$
\left\{
\begin{aligned}
\alpha_2\beta_3 - \alpha_3\beta_2 = 0\\
\alpha_1\beta_3 - \alpha_3\beta_1 = 0\\ 
\alpha_1\beta_2 - \alpha_2\beta_1 = 0\\
\end{aligned}
\right.
$$
hold. If $\alpha_i\neq 0$ for $i\in \{1, 2, 3\}$, then by setting $\beta_i = 0$, we obtain $\alpha_i\beta_j = 0$, where $j\neq i$. Then $\alpha_i = 0$ for $\beta_j = 1$; this is a contradiction. Consequently, $C(\Delta(Q_8, Q'_8)) = \{\alpha_0f \, \mid \, \alpha_0\in C(R)\}$. Therefore, for the element $af\in\Delta(Q_8, Q'_8)\setminus C$, there does not exist an element $c\in C$ such that $0\neq af\cdot c\in C$. Therefore, $\Delta(Q_8, Q'_8)$ is not a centrally essential ring. Therefore, $RQ_8$ is not a centrally essential ring. This is a contradiction.

Conversely, let $R$ be a centrally essential ring of characteristic $2$. It follows from \cite[Theorem 9]{C63} that the fundamental ideal 
$\Delta = \Delta(Q_8, Q_8)$ is nilpotent. Next we use the proof method from the paper \cite[Proposition 2.6]{MT18}. For $0\neq x\in RQ_8$, we consider all possible products $x(1 - z)$, where $z\in Z(Q_8)$. If at least one of them, for example $x_1 = x(1 - z_1)$, is not equal to zero, then we consider $x_1(1 - z)$, and so on. This process will terminate after a finite number of steps, i.e., there exists a positive integer $k$ such that $x_k\cdot\Delta(Q_8, Z(Q_8)) = 0$ (we assume that that $x_0 = x$). Then $x_k\in RQ_8\sum_{Z(Q_8)}$; e.g., see \cite[Proposition 1.12]{ZM75}. Therefore, $RQ_8\sum_{Z(Q_8)}\subseteq C(RQ_8)$. Indeed, if 
$g, h\in Q_8$, then
$$
[g, h\sum_{Z(Q_8)}] = [g, h]\sum_{Z(Q_8)} = gh(1 - h^{-1}g^{-1}hg)\sum_{Z(Q_8)} = 0,
$$
since $Q_8$ is a nilpotent group of nilpotence index $2$ and $h^{-1}g^{-1}hg\in Q'_8 = Z(Q_8)$. By setting $c = (1 - z_1)\ldots (1 - z_k)$ (or $c = 1$, for $k = 0$), we obtain $c\in C(RQ_8)$ and $0\neq xc = x_k\in C(RQ_8)$. Consequently, $RQ_8$ is a centrally essential ring.~$\square$

\textbf{Example 3.7.} Let $R = \mathbb{Z}_9$. We consider the following matrix representation of the ideal $\Delta(Q_8, Q'_8)$:
$$
M_{\Delta} = \begin{pmatrix}
q_0& q_1 & q_2 & q_3\\
- q_1& q_0 & q_3 & - q_2\\
- q_2& - q_3 & q_0 & q_1\\
- q_3& q_ 2& - q_1 & q_0
\end{pmatrix}.
$$
As proven above, $M_{\Delta}$ is not a centrally essential ring. Moreover, $M_{\Delta}$ has non-central idempotents, for example
$$
e = \begin{pmatrix}
5 & 1 & 1 & 0\\
- 1 & 5 & 0 & - 1\\
- 1 & 0 & 5 & 1\\
0 & 1 & - 1 & 5
\end{pmatrix}.
$$

It follows from \cite[Lemma 2.3]{MT18} that $\mathbb{Z}_9Q_8$ is not a centrally essential ring.~$\square$ 

\section{Supplements and Remarks}\label{sec4}

A right module $M$ over a ring $A$ is a \textsf{rational extension} of its submodule $N$ if for every submodule $B$ such that $N\subseteq B\subseteq M$ and for any $f\in Hom_R(B, M)$, the condition $f(N) = 0$ holds if and only if $f = 0$. This is equivalent to the property that for any $x, 0\neq y\in M$, there exist $r\in R$ and $n\in \mathbb{Z}$ such that $xr + xn\in N$ and $yr + yn\neq 0$, e.g., see \cite[Proposition 7.2]{F67}.

A ring $A$ is said to be \textsf{centrally rational} if the module $A_C$ is a rational extension of the module $C_C$. 

\textbf{Remark 4.1.} We note that all centrally rational rings are centrally essential. The converse is not true. Indeed, if $F$ is a field $\mathbb{Z}/2\mathbb{Z}$ and $A$ is a group $F$-algebra of the group quaternion of order 8, then $A$ is a non-commutative centrally essential ring, while every centrally rational ring is commutative, by the following proposition. 

\textbf{Proposition 4.2.} A ring $A$ is commutative if and only if 
$A$ is a centrally rational ring.

\textbf{Proof.} Let $ab\neq ba$ for some $a, b\in R$. Then $ac\in C$ and $(ab - ba)c\neq 0$ for some $c\in C$. Then 
$abc= (ac)b = bac$. This is a contradiction.~$\square$

It is well known that if a ring $R$ is nonsingular, then the ring $Q_{max}(R)$ is injective; e.g., see \cite[Theorem 1+2, P. 69]{F67}. However, if $Z_r(R)\neq 0$, then $Q_{max}(R)$ could also be injective.

\textbf{Example 4.3.} We give an example of a non-commutative centrally essential ring with right injective maximal ring of fractions.

We consider the group algebra $R = \mathbb{Z}_2 Q_8$ with center
$$
C(R) = \{\lambda_0 + \lambda_1a^2 + \lambda_2(a + a^3) + \lambda_3(b + a^3b) + \lambda_4(ab + a^3b),\quad \lambda_i\in \mathbb{Z}_2\};
$$
e.g., see \cite[Part 2]{P77}. Then $|C(R)| = 32$, $R$ is a ring with subdirectly indecomposable center and core $H = \{0, \widehat Q_8\}$, where $\widehat Q_8 = \sum_{g\in Q_8}g$. In addition, $R$ is a centrally essential ring, by Theorem 3.6. It follows from \cite[Corollary 13.5.5(5), Theorem 13.5.7]{K82} that $R$ is a $QF$-ring. Then the module $R_R$ is injective; see \cite[Theorem 13.2.1]{K82}. It follows from \cite[Theorem B, P. 64]{F67} that $Q_{max}(R) = R$. It follows from 4.2 that the ring $R$ is not centrally rational, which can be verified directly. Indeed, let $0\neq a\in H$, $b\in Q_8$. We note that $bc\notin C(R)$ for all invertible central elements $c$ (it is easy to see that there are 16 such elements, and they are sums of an odd number of conjugacy classes as terms). Since $R$ is a centrally essential ring, $bd\in C(R)$ for some non-invertible central element $d$. On the other side, $d$ is a nilpotent element and $dH = 0$; see \cite[Theorem 2.1(v)]{Desh71}. Consequently, $R$ is not a centrally rational ring.~$\square$

The group algebra $\mathbb{Z}_2 Q_8$ contains 256 elements. In the paper \cite[Example 4.1, $\mathbb{F} = \mathbb{Z}_2$]{LT21b}, it is constructed a non-commutative centrally essential ring consisting of 128 elements. At the present time, the authors do not know examples of non-commutative centrally essential rings, containing less than 128 elements. A different situation occurs when considering centrally essential semirings. We recall that a \textsf{semiring} is a structure that differs from an associative ring, because the additive operation is not allways reversible. The zero of a semiring $S$ is multiplicative by definition: for any $s\in S$, we have $0s = s0 = 0$. A semiring $S$ is said to be \textsf{semisubtractive} if for any distinct $a, b\in S$, there exists an element $x\in S$ such that $a + x = b$ or $b + x = a$, e.g., see \cite{HW96}.

\textbf{Example 4.4.} We consider a semiring $S = \{0, 1, a, b, c\}$ in which operations are defined as follows:
$$
\begin{tabular}{ |c|c|c|c|c|c|}
 \hline
$+$ & $0$ & $1$ & $a$ & $b$ & $c$ \\
\hline
$0$ & $0$ & $1$ & $a$ & $b$ & $c$ \\
\hline
$1$ & $1$ & $1$ & $1$ & $b$ & $1$ \\
\hline
$a$ & $a$ & $1$ & $a$ & $b$ & $a$ \\
\hline
$b$ & $b$ & $b$ & $b$ & $b$ & $b$ \\
\hline
$c$ & $c$ & $1$ & $a$ & $b$ & $c$ \\
 \hline
\end{tabular}, \quad \quad
\begin{tabular}{ |c|c|c|c|c|c|}
 \hline
$\cdot$ & $0$ & $1$ & $a$ & $b$ & $c$ \\
\hline
$0$ & $0$ & $0$ & $0$ & $0$ & $0$ \\
\hline
$1$ & $0$ & $1$ & $a$ & $b$ & $c$ \\
\hline
$a$ & $0$ & $a$ & $a$ & $a$ & $c$ \\
\hline
$b$ & $0$ & $b$ & $b$ & $b$ & $c$ \\
\hline
$c$ & $0$ & $c$ & $c$ & $c$ & $c$ \\
 \hline
\end{tabular}.
$$
The center $C(S)$ of the semiring $S$ consists of elements $0$, $1$ and $c$. Since $0\neq x\cdot c = c\in C(S)$ for $0\neq x\in S$, we have that $S$ is a non-commutative semisubtractive centrally essential semiring of order 5.

\end{document}